# HYPOTHESIS TEST FOR NORMAL MIXTURE MODELS: THE EM APPROACH

By Jiahua Chen[1] and Pengfei Li

*University of British Columbia and University of Alberta*

Normal mixture distributions are arguably the most important mixture models, and also the most technically challenging. The likelihood function of the normal mixture model is unbounded based on a set of random samples, unless an artificial bound is placed on its component variance parameter. Moreover, the model is not strongly identifiable so it is hard to differentiate between over dispersion caused by the presence of a mixture and that caused by a large variance, and it has infinite Fisher information with respect to mixing proportions. There has been extensive research on finite normal mixture models, but much of it addresses merely consistency of the point estimation or useful practical procedures, and many results require undesirable restrictions on the parameter space. We show that an EM-test for homogeneity is effective at overcoming many challenges in the context of finite normal mixtures. We find that the limiting distribution of the EM-test is a simple function of the $0.5\chi_0^2 + 0.5\chi_1^2$ and $\chi_1^2$ distributions when the mixing variances are equal but unknown and the $\chi_2^2$ when variances are unequal and unknown. Simulations show that the limiting distributions approximate the finite sample distribution satisfactorily. Two genetic examples are used to illustrate the application of the EM-test.

**1. Introduction.** The class of finite normal mixture models has many applications. More than a hundred years ago, Pearson (1894) modeled a set of crab observations with a two-component normal mixture distribution. In genetics, such models are often used for quantitative traits influenced by major genes. Roeder (1994) discusses an example in which the blood chemical concentration of interest is influenced by a major gene with additive effects [see

Received March 2008; revised August 2008.

[1]Supported in part by the Natural Sciences and Engineering Research Council of Canada and by Mathematics of Information Technology and Complex Systems and the start-up grant of University of Alberta.

Schork, Allison and Thiel (1996)] for additional examples in human genetics. Normal mixture models are also used to account for heterogeneity in the age of onset for male and female schizophrenia patients [Everitt (1996)], and used in hematology studies [McLaren (1996)]. They play a fundamental role in cluster analysis [Tadesse, Sha and Vannucci (2005) and Raftery and Dean (2006)], and in the study of the false discovery rate [Efron (2004), McLachlan, Bean and Ben-Tovim (2006), Sun and Cai (2007) and Cai, Jin and Low (2007)]. In financial economics they are used for daily stock returns [Kon (1984)].

Contrary to intuition, of all the finite mixture models, the normal mixture models have the most undesirable mathematical properties. Their likelihood functions are unbounded unless the component variances are assumed equal or constrained, the Fisher information can be infinite and the strong identifiability condition is not satisfied. We demonstrate these points in the following example.

EXAMPLE 1.  Let $X_1, \ldots, X_n$ be a random sample from the following normal mixture model:

$$(1.1) \qquad\qquad (1-\alpha)N(\theta_1, \sigma_1^2) + \alpha N(\theta_2, \sigma_2^2).$$

Let $f(x, \theta, \sigma)$ be the density function of a normal distribution with mean $\theta$ and variance $\sigma^2$. The likelihood function is given by

$$(1.2) \quad l_n(\alpha, \theta_1, \theta_2, \sigma_1, \sigma_2) = \sum_{i=1}^{n} \log\{(1-\alpha)f(X_i; \theta_1, \sigma_1) + \alpha f(X_i; \theta_2, \sigma_2)\}:$$

1. *Unbounded likelihood function.* The log-likelihood function is unbounded for any given $n$ because when $\theta_1 = X_1$, $0 < \alpha < 1$, it goes to infinity when $\sigma_1 \to 0$.

2. *Infinite Fisher information.* For each given $\theta_1, \theta_2$, $\sigma_1^2$ and $\sigma_2^2$, we have

$$S_n = \frac{\partial l_n(\alpha, \theta_1, \theta_2, \sigma_1, \sigma_2)}{\partial \alpha}\bigg|_{\alpha=0} = \sum_{i=1}^{n} \left\{ \frac{f(X_i; \theta_2, \sigma_2)}{f(X_i; \theta_1, \sigma_1)} - 1 \right\}.$$

   Under the homogeneous model in which $\theta_1 = 0$, $\sigma_1 = 1$ and $\alpha = 0$, that is, $N(0,1)$, the Fisher information

$$E\{S_n^2\} = \infty, \qquad \text{whenever } \sigma_2^2 > 2.$$

3. *Loss of strong identifiability.* It can be seen that

$$\frac{\partial^2 f(x; \theta, \sigma)}{\partial \theta^2}\bigg|_{(\theta, \sigma^2)=(0,1)} = 2\frac{\partial f(x; \theta, \sigma)}{\partial(\sigma^2)}\bigg|_{(\theta, \sigma^2)=(0,1)}.$$

   This is in violation of the strong identifiability condition introduced in Chen (1995).

The above properties of finite normal mixture models are in addition to other undesirable properties of general finite mixture models. In Hartigan (1985), Liu, Pasarica and Shao (2003) and Liu and Shao (2004), the likelihood ratio statistic is shown to diverge to infinity as the sample size increases, which forces most researchers to restrict the mixing parameter ($\theta$) into some compact space. Without which, Hall and Stewart (2005) find the likelihood ratio test can only consistently detect alternative models at distance $(n^{-1} \log \log n)^{1/2}$ rather than at the usual distance $n^{-1/2}$. The partial loss of identifiability, when $\theta_1 = \theta_2$, once forced in a technical separate condition, $|\theta_1 - \theta_2| \geq \varepsilon > 0$ [Ghosh and Sen (1985)]. This condition has recently been shown to be unnecessary by many authors, for instance, Garel (2005).

The unbounded likelihood prevents straightforward use of maximum likelihood estimation. Placing an additional constraint on the parameter space [e.g., Hathaway (1985)] or adding a penalty function [Chen, Tan and Zhang (2008)] to the log likelihood regains the consistency and efficiency of the maximum constrained or penalized likelihood estimators.

The loss of strong identifiability results in a lower best possible rate of convergence [Chen (1995) and Chen and Chen (2003)]. Furthermore, it invalidates many elegant asymptotic results such as those in Dacunha-Castelle and Gassiat (1999), Chen, Chen and Kalbfleisch (2001) and Charnigo and Sun (2004). Finite Fisher information is a common hidden condition of these papers, but it did not gain much attention until the paper of Li, Chen and Marriott (2008).

Due to the indisputable importance of finite normal mixture models, developing valid and useful statistical procedures is an urgent task, particularly for the test of homogeneity. Yet the task is challenging for the reasons presented. Many existing results used simulated quantiles of the corresponding statistics [see Wolfe (1971), McLachlan (1987) and Feng and McCulloch (1994)]. Without rigorous theory, however, it is difficult to reconcile their varying recommendations.

In this paper, we investigate the application of the EM-test [Li, Chen and Marriott (2008)] to finite normal mixture models and show that this test provides a most satisfactory and general solution to the problem. Interestingly, our asymptotic results do not require any constraints on the mean and variance parameters or compactness of the parameter space.

In Section 2, we present the result for the normal mixture model (1.1) when $\sigma_1^2 = \sigma_2^2 = \sigma^2$. The limiting distribution of the EM-test is shown to be a simple function of the $0.5\chi_0^2 + 0.5\chi_1^2$ and the $\chi_1^2$ distributions. In Section 3, we present the result for the general normal mixture model (1.1). The limiting distribution of the EM-test is found to be the $\chi_2^2$. Both results are stunningly simple and convenient to apply. In both cases, we conduct simulation studies and the outcomes are in good agreement with the asymptotic results. In Section 4, we give two genetic examples. For convenience of the presentation, the proofs are outlined in the Appendix and included in a technical report
[Chen and Li (2008)].

## 2. Normal mixture models in the presence of the structural parameter.
When $\sigma_1 = \sigma_2 = \sigma$ and $\sigma$ is unknown in model (1.1), we call $\sigma$ the structural
parameter. We are interested in the test of the homogeneity null hypothesis

$$H_0: \alpha(1-\alpha)(\theta_1 - \theta_2) = 0$$

under this assumption. Without loss of generality, we assume $0 \leq \alpha \leq 0.5$.

Because the population variance $\mathrm{Var}(X_1)$ is the sum of the component
variance $\sigma^2$ and the variance of the mixing distribution $\alpha(1-\alpha)(\theta_1 - \theta_2)^2$,
$\sigma^2$ is often underestimated by straight likelihood methods. Furthermore,
most asymptotic results are obtained by approximating the likelihood func-
tion with some form of quadratic function [Liu and Shao (2003), Marriott
(2007)]. The approximation is most precise when the fitted $\alpha$ value is away
from 0 and 1. Based on these considerations, we recommend using the mod-
ified log likelihood

$$pl_n(\alpha, \theta_1, \theta_2, \sigma) = l_n(\alpha, \theta_1, \theta_2, \sigma, \sigma) + p_n(\sigma) + p(\alpha)$$

with $l_n(\cdot)$ given in (1.2). We usually select $p_n(\sigma)$ such that it is bounded
when $\sigma$ is large, but goes to negative infinity as $\sigma$ goes to 0, and $p(\alpha)$ such
that it is maximized at $\alpha = 0.5$ and goes to negative infinity as $\alpha$ goes to 0
or 1. Concrete recommendations will be given later.

To construct the EM-test, we first choose a set of $\alpha_j \in (0, 0.5], j = 1, 2, \ldots, J$,
and a positive integer $K$. For each $j = 1, 2, \ldots, J$, let $\alpha_j^{(1)} = \alpha_j$ and compute

$$(\theta_{j1}^{(1)}, \theta_{j2}^{(1)}, \sigma_j^{(1)}) = \arg\max_{\theta_1, \theta_2, \sigma} pl_n(\alpha_j^{(1)}, \theta_1, \theta_2, \sigma).$$

For $i = 1, 2, \ldots, n$ and the current $k$, we use an E-step to compute

$$w_{ij}^{(k)} = \frac{\alpha_j^{(k)} f(X_i; \theta_{j2}^{(k)}, \sigma_j^{(k)})}{(1 - \alpha_j^{(k)}) f(X_i; \theta_{j1}^{(k)}, \sigma_j^{(k)}) + \alpha_j^{(k)} f(X_i; \theta_{j2}^{(k)}, \sigma_j^{(k)})}$$

and then update $\alpha$ and other parameters by an M-step such that

$$\alpha_j^{(k+1)} = \arg\max_\alpha \left\{ \left(n - \sum_{i=1}^n w_{ij}^{(k)}\right) \log(1-\alpha) + \sum_{i=1}^n w_{ij}^{(k)} \log(\alpha) + p(\alpha) \right\}$$

and

$$(\theta_{j1}^{(k+1)}, \theta_{j2}^{(k+1)}, \sigma_j^{(k+1)}) = \arg\left[ \max_{\theta_1, \theta_2, \sigma} \sum_{h=1}^2 \sum_{i=1}^n w_{ij}^{(k)} \log\{f(X_i; \theta_h, \sigma)\} + p_n(\sigma) \right].$$

The E-step and the M-step are iterated $K-1$ times.

For each $k$ and $j$, we define

$$M_n^{(k)}(\alpha_j) = 2\{pl_n(\alpha_j^{(k)}, \theta_{j1}^{(k)}, \theta_{j2}^{(k)}, \sigma_j^{(k)}) - pl_n(1/2, \hat{\theta}_0, \hat{\theta}_0, \hat{\sigma}_0)\},$$

where $(\hat{\theta}_0, \hat{\sigma}_0) = \arg\max_{\theta, \sigma} pl_n(1/2, \theta, \theta, \sigma)$.

The EM-test statistic is then defined as

$$EM_n^{(K)} = \max\{M_n^{(K)}(\alpha_j) : j = 1, \ldots, J\}.$$

We reject the null hypothesis when $EM_n^{(K)}$ exceeds some critical value to be determined.

Consider the simplest case where $J = K = 1$ and $\alpha_1 = 0.5$. In this case, the EM-test is the likelihood ratio test against the alternative models with known $\alpha = 0.5$. The removal of one unknown parameter in the model simplifies the asymptotic property of the (modified) likelihood ratio test, and the limiting distribution becomes the $0.5\chi_0^2 + 0.5\chi_1^2$ which does not require the parameter space of $\theta$ to be compact. The price of this simplicity is a loss of efficiency when the data are from an alternative model with $\alpha \neq 0.5$. Choosing $J > 1$ initial values of $\alpha$ reduces the efficiency loss because the true $\alpha$ value can be close to one of the initial values. The EM-iteration updates the value of $\alpha_j$ and moves it toward the true $\alpha$-value while retaining the nice asymptotic property.

Specific choice of initial set of $\alpha$ values is not crucial in general. This is another benefit of the EM-iteration. The updated $\alpha$-values from either $\alpha = 0.3$ or $\alpha = 0.4$ are likely very close after two iterations. Hence, we recommend $\{0.1, 0.3, 0.5\}$. If some prior information indicates that the potential $\alpha$ value under the alternative model is low, then choosing $\{0.01, 0.025, 0.05, 0.1\}$ can improve the power of the test. We do not investigate the potential refinements further but leave them as a future research project at this stage.

The idea of the EM-test was introduced by Li, Chen and Marriott (2008) for mixture models with a single mixing parameter. Yet finite normal mixture models do not fit into the general theory and pose specific technical challenges. The asymptotic properties of the EM-test will be presented in the next section. The recommendation for penalty functions will be given in Section 2.2.

2.1. *Asymptotic properties.* We study the asymptotic properties of the EM-test under the following conditions on the penalty functions $p(\alpha)$ and $p_n(\sigma)$:

C0. $p(\alpha)$ is a continuous function such that it is maximized at $\alpha = 0.5$ and goes to negative infinity as $\alpha$ goes to 0 or 1.

C1. $\sup\{|p_n(\sigma)| : \sigma > 0\} = o(n)$.

C2. The derivative $p_n'(\sigma) = o_p(n^{1/4})$ at any $\sigma > 0$.

We allow $p_n$ to be dependent on the data. To ensure that the EM-test has the invariant property, we recommend choosing a $p_n$ that also satisfies the following:

C3.  $p_n(a\sigma; aX_1 + b, \ldots, aX_n + b) = p_n(\sigma; X_1, \ldots, X_n)$.

The following intermediate results reveal some curious properties of the finite normal mixture model.

THEOREM 1.   *Suppose conditions* C0, C1, *and* C2 *hold. Under the null distribution* $N(\theta_0, \sigma_0^2)$ *we have, for* $j = 1, \ldots, J$ *and any* $k \leq K$, *the following:*

(a)  *if* $\alpha_j = 0.5$, *then:*

$$\theta_{j1}^{(k)} - \theta_0 = O_p(n^{-1/8}), \qquad \theta_{j2}^{(k)} - \theta_0 = O_p(n^{-1/8}),$$

$$\alpha_j^{(k)} - \alpha_j = O_p(n^{-1/4}), \qquad \sigma_j^{(k)} - \sigma_0 = O_p(n^{-1/4}),$$

(b)  *if* $0 < \alpha_j < 0.5$, *then:*

$$\theta_{j1}^{(k)} - \theta_0 = O_p(n^{-1/6}), \qquad \theta_{j2}^{(k)} - \theta_0 = O_p(n^{-1/6}),$$

$$\alpha_j^{(k)} - \alpha_j = O_p(n^{-1/4}), \qquad \sigma_j^{(k)} - \sigma_0 = O_p(n^{-1/3}).$$

Note that the convergence rates of $(\theta_{j1}^{(k)}, \theta_{j2}^{(k)}, \sigma_j^{(k)})$ depend on the choice of initial $\alpha$ value, and it singles out $\alpha = 0.5$. Even when $\alpha_1 = 0.5$, $\alpha_1^{(k)} \neq 0.5$ when $k > 1$. However, this does not reduce case (a) to case (b) because $\alpha_1^{(k)} = 0.5 + o_p(1)$ rather than equaling a nonrandom constant $\alpha_1 \neq 0.5$.

THEOREM 2.   *Suppose conditions* C0, C1, *and* C2 *hold and* $\alpha_1 = 0.5$. *Then, under the null distribution* $N(\theta_0, \sigma_0^2)$ *and for any finite* $K$ *as* $n \to \infty$,

$$\Pr(EM_n^{(K)} \leq x) \to F(x - \Delta)\{0.5 + 0.5F(x)\},$$

*where* $F(x)$ *is the cumulative density function* (*CDF*) *of the* $\chi_1^2$ *and*

$$\Delta = 2 \max_{\alpha_j \neq 0.5} \{p(\alpha_j) - p(0.5)\}.$$

To shed some light on the nonconventional results, we reveal some helpful momental relationships. Without loss of generality, assume that under the null model $\theta_1 = \theta_2 = 0$ and $\sigma^2 = 1$. The EM-test or other likelihood-based methods fit the data from the null model with an alternative model $(1 - \alpha)N(\theta_1, \sigma^2) + \alpha N(\theta_2, \sigma^2)$. Asymptotically, the fit matches the first few sample moments. When $\alpha = 0.5$ is presumed, the first three moments of a homogeneous model and an alternative model can be made identical with proper choice of the values of the remaining parameters. Which model fits the data better is revealed through the fourth moment,

$$E(X_1^4) = 3 - (\theta_1^4 + \theta_2^4) \leq 3.$$

Thus, for local alternatives, we may as well test

$$H_0 : E(X_1^4) = 3 \quad \text{versus} \quad H_a : E(X_1^4) < 3.$$

The parameter of this null hypothesis is on the boundary so that the null limiting distribution of $M_n(0.5)$ is the $0.5\chi_0^2 + 0.5\chi_1^2$.

When $\alpha = \alpha_0 \in (0, 0.5)$, the first two moments of the null and alternative models can be made identical, but their third moments differ because

$$E(X_1^3) = (1 - \alpha_0)\theta_1^3 + \alpha_0\theta_2^3,$$

which can take any value in a neighborhood of 0. Thus, for local alternatives, we may as well test

$$H_0 : E(X_1^3) = 0 \quad \text{versus} \quad H_a : E(X_1^3) \neq 0.$$

Because the null hypothesis is an interior point, $M_n(\alpha_0)$ has the asymptotic distribution $\chi_1^2 + 2\{p(\alpha_0) - p(0.5)\}$ in which $2\{p(\alpha_0) - p(0.5)\}$ is due to the penalty.

Since the sample third and fourth moments are asymptotically orthogonal, the limiting distribution of the EM-test involves the maximum of two independent distributions, the $\chi_1^2$ and the $0.5\chi_0^2 + 0.5\chi_1^2$, and a term caused by the penalty $p(\alpha)$. This is the result as in the above theorem.

The order assessment results in Theorem 1 can be similarly explained. If $\alpha = 0.5$ is presumed, the fitted fourth moment of the mixing distribution will be $O_p(n^{-1/2})$ and hence both fitted $\theta_1$ and $\theta_2$ are $O_p(n^{-1/8})$. For other $\alpha$ values, the fitted third moment is $O_p(n^{-1/2})$, which implies that the fitted $\theta_1$ and $\theta_2$ are $O(n^{-1/6})$.

2.2. *Simulation results.* We demonstrate the precision of the limiting distribution via simulation and explore the power properties. Among several existing results, the modified likelihood ratio test (MLRT) in Chen and Kalbfleisch (2005) is known to have an accurate asymptotic upper bound. Thus, we also include this method in our simulation. The likelihood ratio test (LRT) is included due to its popularity among researchers and simulated its critical values.

The key idea of the MLRT is to define the modified likelihood function as

$$\tilde{l}_n(\alpha, \theta_1, \theta_2, \sigma) = l_n(\alpha, \theta_1, \theta_2, \sigma) + p(\alpha)$$

and the recommended penalty function is $\log\{4\alpha(1-\alpha)\}$. The corresponding statistic is defined as

$$M_n = 2\{l_n(\tilde{\alpha}, \tilde{\theta}_1, \tilde{\theta}_2, \tilde{\sigma}) - l_n(0.5, \tilde{\theta}_0, \tilde{\theta}_0, \tilde{\sigma}_0)\},$$

where $(\tilde{\alpha}, \tilde{\theta}_1, \tilde{\theta}_2, \tilde{\sigma})$ and $(0.5, \tilde{\theta}_0, \tilde{\theta}_0, \tilde{\sigma}_0)$ maximize $\tilde{l}_n$ under the alternative and null models, respectively. Unlike that for the EM-test, the limiting distribution of $M_n$ is unknown but is shown to have an upper bound $\chi_2^2$ when $\theta$ is confined in a compact space. Chen and Kalbfleisch (2005) show that the type I errors of the MLRT with critical values determined by the $\chi_2^2$ distribution are close to the nominal values.

For the EM-test statistics, we choose the penalty function

$$p_n(\sigma) = -\{s_n^2/\sigma^2 + \log(\sigma^2/s_n^2)\},$$

where $s_n^2 = n^{-1}\sum_{i=1}^{n}(X_i - \bar{X})^2$ with $\bar{X} = n^{-1}\sum_{i=1}^{n}X_i$.

It can be seen that (a) $p_n(\sigma)$ satisfies conditions C1–C3; (b) it effectively places an inverse gamma prior on $\sigma^2$; (c) it allows a closed-form expression for $\sigma_j^{(k)}$; and (d) it is maximized at $\sigma^2 = s_n^2$. In fact, even a constant function $p_n(\sigma)$ satisfies C1–C2. This choice of $p_n(\sigma)$ prevents under estimation of $\sigma^2$ and plays a role of higher-order adjustment.

For the penalty function $p(\alpha)$, we choose $p(\alpha) = \log(1 - |1 - 2\alpha|)$. We refer to Li, Chen and Marriott (2008) for reasons of this choice. The combination of $p_n(\sigma)$ and $p(\alpha)$ results in accurate type I errors for the EM-test.

We conducted the simulation with two groups of initial values for $\alpha$: $(0.1, 0.2, 0.3, 0.4, 0.5)$ and $(0.1, 0.3, 0.5)$. We generated 20,000 random samples from $N(0, 1)$ with sample size $n$ ($n = 100, 200$). The simulated null rejection rates are summarized in Table 1. The EM-test and the MLRT both have accurate type I errors, especially $EM_n^{(2)}$ with the three initial values $(0.1, 0.3, 0.5)$ for $\alpha$.

We selected four models for power assessment. The parameter settings are shown in rows 2–5 of Table 2. The powers of the EM-test, the MLRT and the LRT are estimated based on 5,000 repetitions and are presented in Table 3. We used the simulated critical values to ensure fairness of the comparison. The results show that the EM-test statistics based on three initial values have almost the same power as those from five initial values. Combining the type I error results and the power comparison results, we recommend the use of $EM_n^{(2)}$ with three initial values $(0.1, 0.3, 0.5)$ for $\alpha$.

The EM-test has higher power when the mixing proportion $\alpha$ is close to 0.5, while the MLRT statistic performs better when $\alpha$ is close to 0. The powers of the LRT and the MLRT are close under all models. However, the limiting distribution of the EM-test is obtained without any restrictions on the model, while the limiting distribution of the MLRT or of the LRT is unknown, and the upper bound result for the MLRT is obtained under

TABLE 1
*Type* I *errors (%) of the EM-test and the MLRT*

| Level | $EM_n^{(1)}$ | $EM_n^{(2)}$ | $EM_n^{(3)}$ | $EM_n^{(1)}$ | $EM_n^{(2)}$ | $EM_n^{(3)}$ | MLRT |
|---|---|---|---|---|---|---|---|
| | | | | $n = 100$ | | | |
| 10% | 8.9 | 9.1 | 9.2 | 9.2 | 9.9 | 10.2 | 10.9 |
| 5% | 4.6 | 4.8 | 4.8 | 4.6 | 5.1 | 5.3 | 5.7 |
| 1% | 0.9 | 1.0 | 1.0 | 0.9 | 1.0 | 1.1 | 1.2 |
| | | | | $n = 200$ | | | |
| 10% | 9.3 | 9.4 | 9.5 | 9.7 | 10.0 | 10.3 | 9.8 |
| 5% | 4.6 | 4.8 | 4.8 | 4.7 | 5.0 | 5.1 | 5.0 |
| 1% | 1.0 | 1.1 | 1.1 | 0.9 | 1.1 | 1.1 | 1.1 |

Results in columns (2, 3, 4) used $\alpha = (0.1, 0.2, 0.3, 0.4, 0.5)$.
Results in columns (5, 6, 7) used $\alpha = (0.1, 0.3, 0.5)$.

some restrictions. When $\alpha$ is small and some prior information on $\alpha$ value is known, the lower efficiency problem of the EM-test can be easily fixed. We conducted additional simulation by choosing the set of initial $\alpha$-values $\{0.1, 0.05, 0.025, 0.01\}$. In this case, the limiting distribution of the EM-test becomes $\chi_1^2 + 2\{p(0.1) - p(0.5)\}$. When $n = 100$, the power comparison between the EM-test and the MLRT becomes 71.5% versus 71.0% for model III, and 73.1% versus 75% for model IV. Therefore, the EM-test can be refined in many ways to attain higher efficiency. Naturally, a systematic way is preferential and is best left to a future research project.

The other eight models in Table 2 have unequal variances, which are mainly selected for power comparisons in Section 3.3. To examine the im-

TABLE 2
*Parameter values of normal mixture models for power assessment*

| | $1 - \alpha$ | $\theta_1$ | $\theta_2$ | $\sigma_1$ | $\sigma_2$ |
|---|---|---|---|---|---|
| Model I | 0.50 | $-1.15$ | 1.20 | 1.00 | 1.00 |
| Model II | 0.25 | $-1.15$ | 1.15 | 1.00 | 1.00 |
| Model III | 0.10 | $-1.30$ | 1.30 | 1.00 | 1.00 |
| Model IV | 0.05 | $-1.55$ | 1.55 | 1.00 | 1.00 |
| Model V | 0.50 | 0 | 0 | 1.20 | 0.50 |
| Model VI | 0.25 | 0 | 0 | 1.15 | 0.50 |
| Model VII | 0.10 | 0 | 0 | 1.40 | 0.50 |
| Model VIII | 0.05 | 0 | 0 | 1.85 | 0.50 |
| Model IX | 0.50 | 0.75 | $-0.75$ | 1.20 | 0.80 |
| Model X | 0.25 | 0.65 | $-0.65$ | 1.20 | 0.80 |
| Model XI | 0.10 | 0.85 | $-0.85$ | 1.20 | 0.80 |
| Model XII | 0.05 | 1.15 | $-1.15$ | 1.20 | 0.80 |

TABLE 3
*Powers (%) of the EM-test, the MLRT and the LRT at 5% level*

| Model | $EM_n^{(1)}$ | $EM_n^{(2)}$ | $EM_n^{(3)}$ | $EM_n^{(1)}$ | $EM_n^{(2)}$ | $EM_n^{(3)}$ | MLRT | LRT |
|-------|-----------|-----------|-----------|-----------|-----------|-----------|------|-----|
| | | | $n = 100$ | | | | | |
| I   | 53.4 | 53.2 | 52.8 | 53.8 | 53.4 | 53.4 | 45.2 | 45.1 |
| II  | 51.8 | 51.7 | 51.6 | 50.3 | 50.5 | 50.7 | 50.7 | 51.0 |
| III | 51.9 | 52.2 | 52.2 | 50.7 | 51.3 | 51.7 | 59.2 | 59.2 |
| IV  | 49.5 | 51.2 | 51.5 | 50.7 | 51.6 | 52.0 | 63.1 | 64.4 |
| V   | 15.2 | 17.0 | 17.6 | 16.0 | 17.8 | 18.1 | 33.4 | 34.1 |
| IX  | 49.4 | 49.3 | 49.1 | 48.1 | 48.7 | 48.6 | 48.3 | 48.4 |
| | | | $n = 200$ | | | | | |
| I   | 85.2 | 85.2 | 85.1 | 85.3 | 85.4 | 85.3 | 80.1 | 78.6 |
| II  | 85.0 | 84.9 | 84.9 | 84.7 | 84.8 | 84.7 | 84.3 | 83.2 |
| III | 86.0 | 86.1 | 86.1 | 85.7 | 85.8 | 85.9 | 90.9 | 89.8 |
| IV  | 81.4 | 82.3 | 82.5 | 82.5 | 83.1 | 83.2 | 91.1 | 90.2 |
| V   | 23.0 | 25.0 | 25.9 | 24.4 | 26.0 | 26.9 | 52.9 | 55.3 |
| IX  | 82.3 | 82.2 | 82.2 | 81.8 | 82.0 | 82.0 | 82.1 | 81.3 |

Results in columns $(2,3,4)$ used $\alpha = (0.1, 0.2, 0.3, 0.4, 0.5)$.
Results in columns $(5,6,7)$ used $\alpha = (0.1, 0.3, 0.5)$.

portance of the equal variance assumption, we applied the current EM-test designed for finite normal mixture models in the presence of a structural parameter to the data from models V and IX. In some sense, model V is a null model because its two component means are equal, while model IX is an alternative model because its two component means are unequal. It can be seen in Table 3 that the current EM-test has a rightfully low rejection rate against model V. This property is not shared by the MLRT. At the same time, the current EM-test has good power for detecting model IX. In fact, the power is comparable to that of the EM-test designed for finite normal mixture models with unequal variances, to be introduced in the next section. We conclude that when $\sigma_1/\sigma_2$ is close to 1, the power of the current EM-test is not sensitive to the $\sigma_1 = \sigma_2$ assumption.

To explore what happens when $\sigma_1/\sigma_2$ is large, we generated data from model IX with $\sigma_1$ reset to 2.4. The current EM-test rejected the null hypothesis 84% of the time, compared to a 96% rejection rate for the EM-test designed for finite mixture models without an equal variance assumption when $n = 100$. We conclude that when the two component variances are rather different, the current EM-test should not be used. An EM-test designed for finite mixture models without an equal variance assumption is preferred.

### 3. Normal mixture models in both mean and variance parameters.

3.1. *The EM-test procedure.* In this section, we apply the EM-test to the test of homogeneity in the general normal mixture model (1.1) where both $\theta$ and $\sigma$ are mixing parameters. We wish to test

$$H_0 : \alpha(1 - \alpha) = 0 \quad \text{or} \quad (\theta_1, \sigma_1^2) = (\theta_2, \sigma_2^2).$$

Compared to the case where $\sigma$ is a structural parameter, the asymptotic properties of likelihood-based methods become much more challenging because of the unbounded log-likelihood and infinite Fisher information. Especially because of the latter, there exist few asymptotic results for general finite normal mixture models. Interestingly, we find that the EM-test can be directly applied and the asymptotic distribution is particularly simple. However, its derivation is complex.

To avoid the problem of unbounded likelihood, adding a penalty becomes essential in our approach. We define

$$pl_n(\alpha, \theta_1, \theta_2, \sigma_1, \sigma_2) = l_n(\alpha, \theta_1, \theta_2, \sigma_1, \sigma_2) + p_n(\sigma_1) + p_n(\sigma_2) + p(\alpha),$$

where $p_n(\sigma)$, $p(\alpha)$ are the same as before.

The EM-test statistic is constructed similarly. We first choose a set of $\alpha_j \in (0, 0.5], j = 1, 2, \ldots, J$ and a positive integer $K$. For each $j = 1, 2, \ldots, J$, let $\alpha_j^{(1)} = \alpha_j$ and compute

$$(\theta_{j1}^{(1)}, \theta_{j2}^{(1)}, \sigma_{j1}^{(1)}, \sigma_{j2}^{(1)}) = \underset{\theta_1, \theta_2, \sigma_1, \sigma_2}{\arg\max} \, pl_n(\alpha_j^{(1)}, \theta_1, \theta_2, \sigma_1, \sigma_2).$$

For $i = 1, 2, \ldots, n$ and the current $k$, we use the E-step to compute

$$w_{ij}^{(k)} = \frac{\alpha_j^{(k)} f(X_i; \theta_{j2}^{(k)}, \sigma_j^{(k)})}{(1 - \alpha_j^{(k)}) f(X_i; \theta_{j1}^{(k)}, \sigma_j^{(k)}) + \alpha_j^{(k)} f(X_i; \theta_{j2}^{(k)}, \sigma_j^{(k)})}$$

and then we use the M-step to update $\alpha$ and other parameters such that

$$\alpha_j^{(k+1)} = \underset{\alpha}{\arg\max} \left\{ \left( n - \sum_{i=1}^{n} w_{ij}^{(k)} \right) \log(1 - \alpha) + \sum_{i=1}^{n} w_{ij}^{(k)} \log(\alpha) + p(\alpha) \right\}$$

and

$$(\theta_{j1}^{(k+1)}, \theta_{j2}^{(k+1)}, \sigma_{j1}^{(k+1)}, \sigma_{j2}^{(k+1)})$$
$$= \underset{\theta_1, \theta_2, \sigma_1, \sigma_2}{\arg\max} \sum_{h=1}^{2} \left[ \sum_{i=1}^{n} w_{ij}^{(k)} \log\{f(X_i; \theta_h, \sigma_h)\} + p_n(\sigma_h) \right].$$

The E-step and the M-step are iterated $K - 1$ times.

For each $k$ and $j$, we define

$$M_n^{(k)}(\alpha_j) = 2\{pl_n(\alpha_j^{(k)}, \theta_{j1}^{(k)}, \theta_{j2}^{(k)}, \sigma_{j1}^{(k)}, \sigma_{j2}^{(k)}) - pl_n(1/2, \hat{\theta}_0, \hat{\theta}_0, \hat{\sigma}_0, \hat{\sigma}_0)\},$$

where $(\hat{\theta}_0, \hat{\sigma}_0) = \arg\max_{\theta,\sigma} pl_n(1/2, \theta, \theta, \sigma, \sigma)$. The EM-test statistic is then defined as

$$EM_n^{(K)} = \max\{M_n^{(K)}(\alpha_j) : j = 1, \ldots, J\}.$$

We reject the null hypothesis when $EM_n^{(K)}$ exceeds some critical value to be determined.

In terms of statistical procedure, the EM-test for the case of $\sigma_1^2 = \sigma_2^2$ is a special case of $\sigma_1^2 \neq \sigma_2^2$. However, the asymptotic distributions and their derivations are different.

3.2. *Asymptotic properties.* We further require that $p_n(\sigma)$ satisfies C1 and:

C4. $p_n'(\sigma) = o_p(n^{1/6})$, for all $\sigma > 0$.
C5. $p_n(\sigma) \leq 4(\log n)^2 \log(\sigma)$, when $\sigma \leq n^{-1}$ and $n$ is large.

The following theorems consider the consistency of $(\alpha_j^{(k)}, \theta_{j1}^{(k)}, \theta_{j2}^{(k)}, \sigma_{j1}^{(k)}, \sigma_{j2}^{(k)})$ and give the major result.

THEOREM 3. *Suppose conditions* C0, C1 *and* C4–C5 *hold. Under the null distribution* $N(\theta_0, \sigma_0^2)$ *we have, for* $j = 1, \ldots, J$, $h = 1, 2$ *and any* $k \leq K$,

$$\alpha_j^{(k)} - \alpha_j = o_p(1), \qquad \theta_{jh}^{(k)} - \theta_0 = o_p(1) \quad and \quad \sigma_{jh}^{(k)} - \sigma_0 = o_p(1).$$

THEOREM 4. *Suppose conditions* C0, C1 *and* C4–C5 *hold. When* $\alpha_1 = 0.5$, *under the null distribution* $N(\theta_0, \sigma_0^2)$ *and for any finite* $K$ *as* $n \to \infty$,

$$EM_n^{(K)} \xrightarrow{d} \chi_2^2.$$

It is a surprise that the EM-test has a simpler limiting distribution when applied to a more complex model. We again shed some light on this via some moment consideration.

The test of homogeneity is to compare the fit of the null $N(0, 1)$ and the fit of the full model. The limiting distribution amounts to considering this problem when the data are from the null model. By matching the first two moments of the full model to the first two sample moments, we roughly select a full model such that

$$(1 - \alpha)\theta_1 + \alpha\theta_2 = 0 \quad \text{and} \quad (1 - \alpha)(\theta_1^2 + \sigma_1^2) + \alpha(\theta_2^2 + \sigma_2^2) = 1.$$

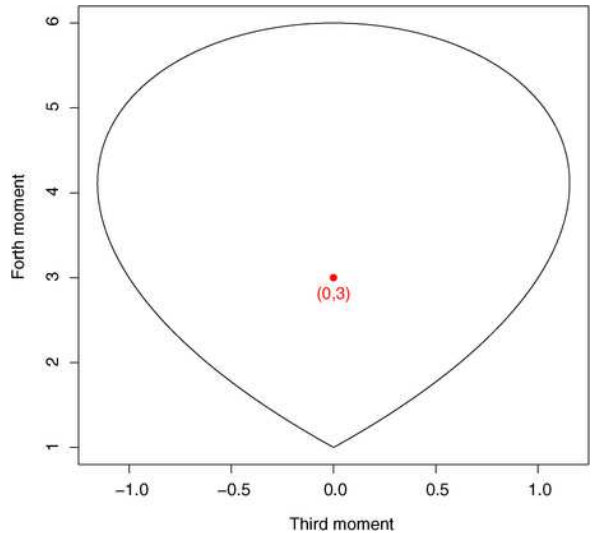

Fig. 1. *The range (area inside the solid line) of $\{E(X_1^3), E(X_1^4)\}$.*

Let $\beta_1 = \theta_1^2 + \sigma_1^2 - 1$. When the value of $\alpha = \alpha_0 \in (0, 0.5]$ (say $\alpha_0 = 0.5$), the third moment and the fourth moment of the full model are

$$E(X_1^3) = 3\theta_1\beta_1,$$
$$E(X_1^4) = 3\beta_1^2 - 2\theta_1^4 + 3.$$

It is easy to verify that $\{E(X_1^3), E(X_1^4)\} = \{0, 3\}$ if and only if the mixture model is the homogeneous model. Therefore, we may as well test

$$H_0 : \{E(X_1^3), E(X_1^4)\} = \{0, 3\} \quad \text{versus} \quad H_a : \{E(X_1^3), E(X_1^4)\} \neq \{0, 3\}.$$

As shown in Figure 1, $\{0, 3\}$ is an interior point of the parameter space of $\{E(X_1^3), E(X_1^4)\}$. Therefore, the null limiting distribution of the EM-test is the $\chi_2^2$. We note that when the observations are from an alternative model, the situation is totally different. A test on moments is not equivalent to the EM-test.

3.3. *Simulation studies.* We demonstrate the precision of the limiting distribution and explore the power properties via simulations. In contrast to the case where $\sigma_1^2 = \sigma_2^2$, the EM-test does not have many competitors. Thus, we set up an MLRT method with

$$M_n = 2\left\{\sup_{\alpha, \theta_1, \theta_2, \sigma_1, \sigma_2} pl_n(\alpha, \theta_1, \theta_2, \sigma_1, \sigma_2) - pl_n(0.5, \hat{\theta}_0, \hat{\theta}_0, \hat{\sigma}_0, \hat{\sigma}_0)\right\}.$$

Although the limiting distribution of $M_n$ is not available, we simulate the critical values and use the MLRT as an efficiency barometer.

We suggest using the penalty function $p_n(\sigma) = -0.25\{s_n^2/\sigma^2 + \log(\sigma^2/s_n^2)\}$, which is almost the same as before except for the coefficient because we have two penalty terms in this problem. Our simulation shows that this choice works well in terms of providing accurate type I errors. We use $p(\alpha) = \log(1 - |1 - 2\alpha|)$ according to the recommendation of Li, Chen and Marriott (2008).

In the simulations, the type I errors were calculated based on 20,000 samples from $N(0, 1)$. As in Section 2.2, we used two groups of initial values $(0.1, 0.2, 0.3, 0.4, 0.5)$ and $(0.1, 0.3, 0.5)$ to calculate $EM_n^{(K)}$. The simulation results are summarized in Table 4. The EM-test statistics based on $(0.1, 0.3, 0.5)$ give accurate type I errors.

The powers of the EM-test and the MLRT for the models in Table 2 are calculated based on 5,000 repetitions and presented in Table 5. Since the limiting distribution of the MLRT is unavailable and hence is not a viable method, the simulated critical values were used for power calculation. The simulation results show that the $EM_n^{(2)}$ and $EM_n^{(3)}$ based on three initial values $(0.1, 0.3, 0.5)$ for $\alpha$ have almost the same power as the MLRT. Further increasing the number of iterations or the number of initial values for $\alpha$ does not increase the power of the EM-test statistics. We therefore recommend the use of $EM_n^{(2)}$ or $EM_n^{(3)}$ based on three initial values $(0.1, 0.3, 0.5)$ for $\alpha$.

We note that when $\sigma_1 = \sigma_2$, the current EM-test loses some power compared to the EM-test designed for finite mixture models in the presence of a structural parameter if the mixing parameter $\alpha$ is close to 0.5, but it has higher power when $\alpha$ is near 0 or 1. Nevertheless, we recommend the use of the current EM-test if the equal variance assumption is likely violated.

## 4. Genetic applications.

TABLE 4
*Type* I *errors (%) of the EM-test*

| Level | $EM_n^{(1)}$ | $EM_n^{(2)}$ | $EM_n^{(3)}$ | $EM_n^{(1)}$ | $EM_n^{(2)}$ | $EM_n^{(3)}$ |
|---|---|---|---|---|---|---|
| | | | $n = 100$ | | | |
| 10% | 10.8 | 10.9 | 10.9 | 10.5 | 10.6 | 10.6 |
| 5% | 5.5 | 5.5 | 5.6 | 5.3 | 5.4 | 5.4 |
| 1% | 1.2 | 1.2 | 1.2 | 1.1 | 1.2 | 1.2 |
| | | | $n = 200$ | | | |
| 10% | 10.7 | 10.7 | 10.7 | 10.4 | 10.5 | 10.5 |
| 5% | 5.4 | 5.4 | 5.4 | 5.1 | 5.2 | 5.2 |
| 1% | 1.1 | 1.1 | 1.1 | 1.0 | 1.0 | 1.0 |

Results in columns $(2, 3, 4)$ used $\alpha = (0.1, 0.2, 0.3, 0.4, 0.5)$.
Results in columns $(5, 6, 7)$ used $\alpha = (0.1, 0.3, 0.5)$.

TABLE 5
*Powers (%) of the EM-test and the MLRT at the 5% level*

| Model | $EM_n^{(1)}$ | $EM_n^{(2)}$ | $EM_n^{(3)}$ | $EM_n^{(1)}$ | $EM_n^{(2)}$ | $EM_n^{(3)}$ | MLRT |
|---|---|---|---|---|---|---|---|
| | | | | $n = 100$ | | | |
| I | 44.0 | 44.0 | 43.9 | 44.1 | 43.8 | 43.8 | 44.0 |
| II | 47.7 | 47.9 | 47.8 | 47.5 | 47.5 | 47.4 | 47.9 |
| III | 55.5 | 55.5 | 55.4 | 55.6 | 55.6 | 55.5 | 55.5 |
| IV | 56.9 | 56.9 | 56.8 | 57.4 | 56.9 | 56.8 | 56.8 |
| V | 58.6 | 58.4 | 58.4 | 58.8 | 58.8 | 58.7 | 58.2 |
| VI | 63.5 | 63.3 | 63.3 | 63.7 | 63.6 | 63.6 | 63.2 |
| VII | 66.8 | 66.6 | 66.6 | 66.9 | 66.8 | 66.8 | 66.6 |
| VIII | 67.3 | 67.2 | 67.2 | 67.4 | 67.2 | 67.1 | 67.4 |
| IX | 48.9 | 48.8 | 48.7 | 49.1 | 48.8 | 48.7 | 48.7 |
| X | 54.6 | 54.6 | 54.5 | 55.0 | 54.8 | 54.6 | 54.3 |
| XI | 56.5 | 56.5 | 56.5 | 57.0 | 56.6 | 56.6 | 56.3 |
| XII | 57.1 | 57.1 | 57.0 | 57.4 | 57.1 | 57.0 | 57.0 |
| | | | | $n = 200$ | | | |
| I | 78.3 | 78.2 | 78.2 | 78.3 | 78.2 | 78.2 | 78.2 |
| II | 82.0 | 81.9 | 81.9 | 82.2 | 82.1 | 82.1 | 81.9 |
| III | 88.6 | 88.6 | 88.6 | 88.8 | 88.7 | 88.7 | 88.5 |
| IV | 88.7 | 88.6 | 88.6 | 88.9 | 88.8 | 88.8 | 88.5 |
| V | 90.0 | 89.9 | 89.9 | 90.1 | 90.0 | 90.0 | 89.8 |
| VI | 91.6 | 91.5 | 91.5 | 91.7 | 91.6 | 91.6 | 91.5 |
| VII | 91.4 | 91.3 | 91.3 | 91.5 | 91.4 | 91.4 | 91.3 |
| VIII | 89.6 | 89.6 | 89.6 | 89.6 | 89.5 | 89.5 | 89.7 |
| IX | 81.7 | 81.5 | 81.5 | 81.9 | 81.8 | 81.7 | 81.4 |
| X | 88.1 | 88.0 | 88.0 | 88.3 | 88.2 | 88.1 | 87.9 |
| XI | 86.4 | 86.3 | 86.3 | 86.5 | 86.4 | 86.4 | 86.2 |
| XII | 87.7 | 87.6 | 87.6 | 87.9 | 87.9 | 87.8 | 87.5 |

Results in columns $(2,3,4)$ used $\alpha = (0.1, 0.2, 0.3, 0.4, 0.5)$.
Results in columns $(5,6,7)$ used $\alpha = (0.1, 0.3, 0.5)$.

EXAMPLE 2. We apply the EM-test to the example discussed in Loisel et al. (1994). Due to the potential use for hybrid production, cytoplasmic male sterility in plant species is a trait of much scientific and economic interest. To efficiently use this character, it is important to find nuclear genes—preferably dominant ones—that induce fertility restoration [MacKenzie and Bassett (1987)]. Loisel et al. (1994) carried out an experiment for detecting a major restoration gene. In this experiment, 150 F2 bean plants were obtained. The number of pods with one up to a maximum of ten grains were then counted on each F2 plant. Loisel et al. (1994) suggested analyzing the square root of the total number of grains for each plant. If a major restoration gene exists, the normal mixture model will provide a more suitable fit; otherwise, the single normal distribution best fits the data. The histogram of the trans-formed counts is given in Figure 2. It indicates the existence of two modes, and an unequal variance normal mixture model is a good choice.

Based on some genetic background, Loisel et al. (1994) postulated a three-component normal mixture model

$$\frac{1}{4}N(\theta_1, \sigma^2) + \frac{1}{2}N(\theta_2, \sigma^2) + \frac{1}{4}N(\theta_3, \sigma^2) \tag{4.1}$$

and tested the null hypothesis that $\theta_1 = \theta_2 = \theta_3$. They found that the limiting distribution of the LRT statistic is a 50–50 mixture of the $\chi_1^2$ and $\chi_2^2$, and the resulting $p$-value is 0.002%. We investigated the null rejection rates of the LRT under model (4.1) when $n = 150$ and the critical values were determined by a 50–50 mixture of the $\chi_1^2$ and $\chi_2^2$ limiting distributions. Based on 40,000 repetitions, the simulated null rejection rates were 15.6%, 8.8% and 2.2% for nominal values of 10%, 5% and 1%. The above $p$-value may be biased toward the liberal side.

For illustration purposes, we re-analyzed the data with the EM-test under model (1.1) with $\sigma_1^2 = \sigma_2^2$. The $p$-value of the MLRT calibrated with the $\chi_2^2$ distribution was found to be 1.4%. We found $EM_n^{(2)} = 6.827$ with three initial values $(0.1, 0.3, 0.5)$ for $\alpha$, corresponding to the $p$-value 1.0%. It can be seen that the EM-test provides stronger evidence against the null model than the MLRT test.

It appears that the equal variance assumption is not suitable. We consider the EM-test for a finite normal mixture with unequal variance. We found that $EM_n^{(1)} = 15.966$ and $EM_n^{(2)} = 20.590$ with three initial values $(0.1, 0.3, 0.5)$ for $\alpha$, resulting in the $p$-values 0.03% and 0.003%, respectively. Further iteration does not change the $p$-value much. This result is in line with the outcome of Loisel et al. (1994). The modified MLES of $(\alpha, \theta_1, \theta_2, \sigma_1, \sigma_2)$ are $(0.175, 10.663, 2.535, 3.203, 1.080)$, confirming that $\sigma_1 \neq \sigma_2$ and explaining why the EM-test under the general model gives much stronger evidence against the null model.

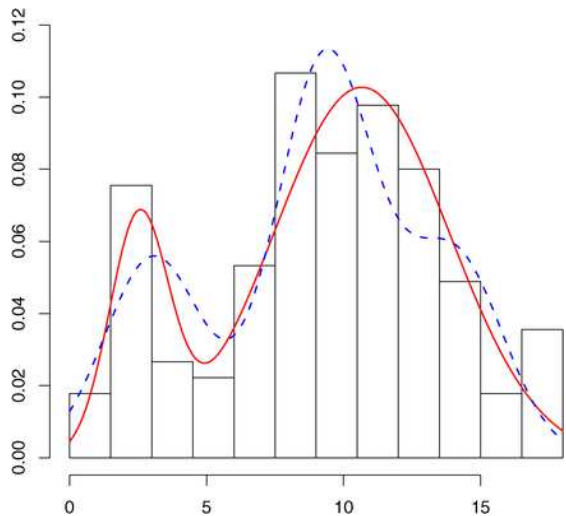

Fig. 2. *The histogram of the square root of the total number of grains per plant, the fitted densities of the normal mixtures in (1.1) (solid line) and in (4.1) (dashed line).*

Figure 2 shows the fitted density functions of models (1.1) and (4.1). Our analysis indicates that a two-component mixture model can fit the data just as well as the model suggested by Loisel et al. (1994). The question of which model is more appropriate is not the focus of this paper.

EXAMPLE 3. The second example considers the data presented in Everitt, Landau and Leese (2001); see part (b) of Table 6.2. This data set is from a schizophrenia study reported by Levine (1981), who collated the results of seven studies on the age of onset of schizophrenia including 99 females and 152 males. We use the male data to illustrate the use of the EM-test. As suggested by Levine (1981), there are two types of schizophrenia in males. The first type is diagnosed at a younger age and is generally more severe; the second type is diagnosed later in life. We wish to test the existence of the two types of schizophrenia.

Everitt, Landau and Leese (2001) fitted the 152 observations using a two-component normal mixture model, and used the LRT to test the homogeneity. Using the $\chi_3^2$ distribution for calibration, they found the $p$-value was less than 0.01%. Following Everitt (1996), our analysis is based on logarithmic transformed data. Assuming model (1.1) with $\sigma_1^2 = \sigma_2^2$, the $p$-value of the MLRT calibrated with the $\chi_2^2$ distribution is 1.8%, but $EM_n^{(2)} = 0$ with three initial values $(0.1, 0.3, 0.5)$ for $\alpha$.

Removing the $\sigma_1 = \sigma_2$ assumption, we find that $EM_n^{(1)} = 13.301$ and $EM_n^{(2)} = 13.323$ with three $\alpha$ initial values $(0.1, 0.3, 0.5)$ and both $p$-values are 0.1%. The modified MLEs of $(\alpha, \mu_1, \mu_2, \sigma_1, \sigma_2)$ are $(0.448, 1.379, 1.319, 0.192, 0.071)$. Our analysis indicates that there are two subpopulations in the population with close mean ages of onset but different variances. This also explains why the EM-test designed for finite mixture models in the presence of a structural parameter is insignificant.

Figure 3 contains the histogram and the fitted densities. It can be seen that the mixture model with unequal variances fits better. We also computed the LRT statistic which equals 15.27 under the unequal variance assumption. If it is calibrated with the $\chi_4^2$ distribution, as suggested by Wolfe (1971), the $p$-value is 0.4%, and if calibrated with the $\chi_6^2$ distribution, as suggested by McLachlan (1987), the $p$-value is 1.8%. Without a solid theory, it would be hard to reconcile these inconsistent outcomes.

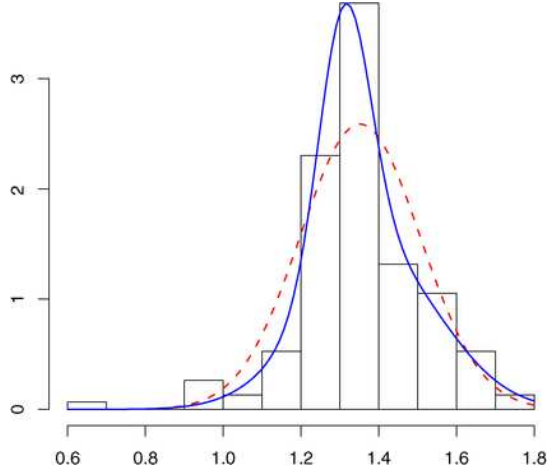

Fig. 3.  *The histogram of log age of onset for male schizophrenics, the fitted densities of the single normal model (dashed line) and normal mixture model in (1.1) (solid line).*

## APPENDIX: PROOFS

To save space, detailed proofs will be given in a technical report [Chen and Li (2008)]. We only include an outline in this appendix.

For both limiting distribution conclusions in Theorems 2 and 4, the first key step is to show that any estimator with $\alpha$ bounded away from 0 and 1, and with a large likelihood value, is consistent for $\theta$ and $\sigma$ under the null model, whether under the equal variance or unequal variance assumption.

The second key step is to show that after finite number of EM-iterations, $\alpha^{(k)}$ keeps a nondeminishing distance from 0 and 1, and the likelihood value is large. Hence, the conclusion in the first key step is applicable.

With both fitted values $(\theta_{1j}^{(k)}, \theta_{2j}^{(k)})$ in a small neighborhood of $\theta_0$, the true value under the null model, the likelihood function is approximated by a quadratic function. This is the third key step.

When $\alpha_j = 0.5$, under the equal variance model, the quadratic approximation leads to the expansion

$$M_n^{(K)}(0.5) = \frac{\{(\sum_{i=1}^n V_i)^-\}^2}{\sum_{i=1}^n V_i^2} + o_p(1).$$

When $\alpha_j \neq 0.5$, the expansion becomes

$$M_n^{(K)}(\alpha_j) = \frac{(\sum_{i=1}^n U_i)^2}{\sum_{i=1}^n U_i^2} + 2\{p(\alpha_j) - p(0.5)\} + o_p(1).$$

Therefore,

$$EM_n^{(K)} = \max\{M_n^{(K)}(\alpha_j), j = 1, 2, \ldots, J\}$$

$$= \max\left[\frac{(\sum_{i=1}^n U_i)^2}{\sum_{i=1}^n U_i^2} + \Delta, \frac{\{(\sum_{i=1}^n V_i)^-\}^2}{\sum_{i=1}^n V_i^2}\right] + o_p(1).$$

We omit the definitions of $U_i$ and $V_i$ but point out that $\sum_{i=1}^n U_i/\sqrt{n}$ and $\sum_{i=1}^n V_i/\sqrt{n}$ are jointly asymptotical bivariate normal and independent. Consequently, the limiting distribution of $EM^{(K)}$ is given by $F(x - \Delta)\{0.5 + 0.5F(x)\}$ with $F(x)$ being the CDF of the $\chi_1^2$ distribution.

Under the unequal variance assumption, the asymptotic expansion is found to be

$$EM_n^{(K)} = \frac{(\sum_{i=1}^n U_i)^2}{\sum_{i=1}^n U_i^2} + \frac{(\sum_{i=1}^n V_i)^2}{\sum_{i=1}^n V_i^2} + o_p(1).$$

Consequently, the limiting distribution of $EM_n^{(K)}$ is the $\chi_2^2$.

Department of Statistics
University of British Columbia
Vancouver
British Columbia, V6T 1Z2
Canada
E-mail: jhchen@stat.ubc.ca

Department of Mathematical
    and Statistical Sciences
University of Alberta
Edmonton
T6G 2G1
Canada
E-mail: pengfei@stat.ualberta.ca